\documentclass[10pt]{article}
\usepackage[latin1]{inputenc}
\usepackage[english]{babel}
\usepackage{amssymb,amsfonts, amsmath}
\usepackage{amsthm}

\usepackage{a4wide}
\newtheorem{theorem}{Theorem}[section]
\newtheorem{proposition}[theorem]{Proposition}
\newtheorem{corollary}[theorem]{Corollary}
\newtheorem{lemma}[theorem]{Lemma}
\newtheorem{remark}[theorem]{Remark}

\newcommand{\bcl}{\begin{center}}
\newcommand{\ecl}{\end{center}}
\newcommand{\brl}{\begin{right}}
\newcommand{\erl}{\end{right}}
\newcommand{\ben}{\begin{enumerate}}
\newcommand{\een}{\end{enumerate}}
\newcommand{\barr}{\begin{array}}
\newcommand{\earr}{\end{array}}
\newcommand{\btab}{\begin{tabular}}
\newcommand{\etab}{\end{tabular}}
\newcommand{\bdoc}{\begin{document}}
\newcommand{\edoc}{\end{document}}
\newcommand{\beqy}{\begin{eqnarray}}
\newcommand{\eeqy}{\end{eqnarray}}

\newcommand{\beqi}{\begin{eqnarray*}}
\newcommand{\eeqi}{\end{eqnarray*}}
\newcommand{\bitem}{\begin{itemize}}

\newcommand{\eitem}{\end{itemize}}
\newcommand{\nln}{\newline}
\newcommand{\newt}{\newtheorem}
\newcommand{\pa}{\partial}
\newcommand{\re}{{I\!\!R}}
\newcommand{\ren}{\re^N}
\newcommand{\xr}{x\in\re }
\newcommand{\x}{\times}
\newcommand{\dyle}{\displaystyle}
\newcommand{\ene}{{I\!\!N}}
\newcommand{\irn}{\int\limits_{\re^N}}
\newcommand{\io}{\int\limits_{\O}}
\newcommand{\meas}{{\rm meas\,}}

\newcommand{\sign}{{\rm sign}}
\newcommand{\map}{\longrightarrow }
\newcommand{\imp}{\Longrightarrow }
\renewcommand{\div}{\nabla\cdot }
\newcommand{\sen}{{\rm sen\,}}
\newcommand{\tg}{{\rm tg\,}}
\newcommand{\arcsen}{{\rm arcsen\,}}
\newcommand{\arctg}{{\rm arctg\,}}
\newcommand{\supp}{{\textsl supp\ }}
\newcommand{\ity}{\int_{-\iy}^{+\iy}}
\newcommand{\limit}{\lim\limits}
\newcommand{\limi}{\limit_{n\to\infty}}
\newcommand{\sumi}{\sum\limits_{n=1}^{\infty}}
\newcommand{\ulu}{\underline u}
\newcommand{\ulw}{\underline w}
\newcommand{\ulz}{\underline z}
\newcommand{\ulv}{\underline v}
\newcommand{\uls}{\underline s}
\newcommand{\olu}{\overline u}
\newcommand{\olv}{\overline v}
\newcommand{\ols}{\overline s}
\newcommand{\ob}{\overline\b}
\newcommand{\ovar}{\overline\var}
\newcommand{\wv}{\widetilde v}
\newcommand{\wu}{\widetilde u}
\newcommand{\ws}{\widetilde s}
\renewcommand{\a }{\alpha }
\renewcommand{\b }{\beta }
\newcommand{\g }{\gamma}
\newcommand{\G }{\Gamma }
\renewcommand{\d }{\delta }

\newcommand{\D }{\Delta }
\newcommand{\e }{\varepsilon }
\newcommand{\z }{\zeta }
\renewcommand{\l }{\lambda }
\renewcommand{\L }{\Lambda }
\newcommand{\m }{\mu }
\newcommand{\n }{\nabla }
\newcommand{\s }{\sigma }
\newcommand{\Sig }{\Sigma }
\renewcommand{\t }{\tau }
\newcommand{\var }{\varphi }
\renewcommand{\o }{\omega }
\renewcommand{\O }{\Omega }
\newcommand{\bR}{{\bf R}}
\newcommand{\bC}{{\bf C}}
\newcommand{\bZ}{{\bf Z}}
\newcommand{\bN}{{\bf N}}
\newcommand{\bQ}{{\bf Q}}
\newcommand{\bK}{{\bf K}}
\newcommand{\bI}{{\bf I}}
\newcommand{\bv}{{\bf v}}
\newcommand{\bV}{{\bf V}}

\newcommand{\pair}[1]{g\left(#1\right)}

\usepackage{pdfsync}

\def\qed{\unskip\kern 6pt \penalty 500
\raise -2pt\hbox{\vrule \vbox to10pt{\hrule width 4pt
\vfill\hrule}\vrule}\par}

\newenvironment{Proof}{\removelastskip\vskip12pt
plus 1pt \noindent\em\rm}{\hfill {\qed \hskip .2cm}}

\title
{Global Solutions of Semilinear \\ Parabolic Equations with Drift Term\\  on Riemannian Manifolds}
\author{Fabio Punzo\thanks{Dipartimento di Matematica, Politecnico di Milano, Italia (fabio.punzo@polimi.it).}}
\date{}

\begin{document}
 \maketitle
{\abstract{ \noindent  We study existence and non-existence of global solutions to the semilinear heat equation with a drift term and a power-like source term $u^p$, on Cartan-Hadamard  manifolds. Under suitable assumptions on Ricci and sectional curvatures, we show that, for any $p>1$, global solutions cannot exists if the initial datum is large enough. Furthermore, under appropriate conditions on the drift term, global existence is obtained for any $p>1$, if the initial datum is sufficiently small. We also deal with Riemannian manifolds whose Ricci curvature tends to zero at infinity sufficiently fast. We show that for any non trivial initial datum, for certain $p$ depending on the Ricci curvature's bound, global solutions cannot exist. On the other hand, for determined values of $p$, depending on the vector field $b$, global solutions exist, for sufficiently small initial data.  
\bigskip


}}

\bigskip
\medskip
\smallskip

\section{Introduction} \setcounter{equation}{0}
We investigate existence and non-existence of nonnegative global 
solutions of Cauchy problems for semilinear parabolic equations of the following
form:
\begin{equation}\label{e1}
\left\{
\begin{array}{ll}
 \,  \pa_t u = \Delta u\, +\, \langle b(x), \nabla u\rangle + u^p \, &\textrm{in}\,\,M\times (0,T)
\\&\\
\textrm{ }u \, = u_0& \textrm{in}\,\,  M\times \{0\} \,;
\end{array}
\right.
\end{equation}
here $M$ is a 
Cartan-Hadamard manifold of dimension $n$, endowed with a metric tensor $g$, $\Delta$ and $\nabla$ denote the Laplace-Beltrami operator and, respectively, the gradient with respect to $g$, $p>1, b$ is a vector field defined on $M$, $u_0$ is a given nonnegative bounded initial datum. For all vector fields $X, Y$ belonging to the tangent bundle of $M$, we set
\[\langle X, Y\rangle := g(X, Y)\,.\] Here $T>0$ is the maximal existence time; when $T=+\infty$, we say that the solution is global. Problem \eqref{e1} can be regarded as a model for a reaction-diffusion process, which takes place on a Riemannian manifold, where $u$ stands for the temperature, $u^p$ is a nonlinear source term and the drift $\langle b, \nabla u\rangle$ can be determined by an external flow field.

\smallskip

Problem \eqref{e1} with $b\equiv 0$ has been widely studied in the literature, both in the Euclidean space  (see, e.g. \cite{Fuji}, \cite{Hay}, \cite{KST}, \cite{Lev}) and on Riemannian manifolds (see, e.g., \cite{GMuP},  \cite{MMP}, \cite{GMP1}, \cite{GMP2}, \cite{Pu12}, \cite{Pu21}, \cite{WY}, \cite{Zhang}). In particular, in \cite{BPT} the hyperbolic space is considered. Some results of \cite{BPT} have been generalized to  Cartan-Hadamard manifolds with strictly negative sectional curvature in \cite{Pu12}, \cite{Pu21}. It is showed that for any $p>1$, if $u_0$ is small enough, then problem \eqref{e1} with $b \equiv 0$ admits a global nonnegative bounded solution. 

Problem \eqref{e1} with $M=\mathbb R^n$ has been investigated, e.g.,  in \cite{AE, BL}. In particular, in \cite{BL} it is proved that for any $p>1$, for sufficiently large initial conditions, global solutions cannot exist. %Furthermore, when, for some $\gamma>-n$ and $r_0>0$,
%$$b(x)=\nabla \phi(x), \,\, \nabla \phi(x)\leq \gamma \quad \text{for all } x\in B_{r_0}^c,$$
%then, for any $1<p<1+\frac{2}{n+\gamma}$, the solution blows-up in finite time (see \cite[Theorem 3.3.1]{BL}. 
On the other hand, if, for some $\nu>-n$,
 \[\langle b(x), x \rangle \geq \nu\quad \text{for all } x\in \mathbb R^n,\]
 and 
 \[p> 1+ \frac 2{n+\nu},\]
then problem \eqref{e1} has a global solution (see \cite[Theorem 3.3.3]{BL}. In \cite{BL} also more general operators are considered and many other results are established.

\medskip

In the present paper, we consider: $(j)$ Cartan-Hadamard manifolds with radial Ricci curvature bounded from below and strictly negative sectional curvature; $(jj)$ Cartan-Hadamard manifolds whose radial Ricci curvature can be negative, but tends to $0$ at infinity fast enough.    

Consider the case $(j)$. We prove that, under suitable assumptions on $b$, if the initial data $u_0$ is big enough, then for any $p>1$ global solutions to problem \eqref{e1} cannot exist.  To prove this result, we use a modification of the method exploited in \cite{Kap} (see also \cite{BL}). 
Note that the curvature assumptions give certain bounds on the volume's growth of geodesics balls, which are essential to employ this method on a complete non-compact Riemannian manifold. Moreover, we will construct and use a function $\varphi\in C^2(M\setminus \partial B_{R_0})\cap C^1(M)$ satisfying 
\begin{equation}\label{e4}
\Delta \varphi - \langle b(x), \nabla \varphi\rangle + \big[\lambda - \operatorname{div} b(x)\big] \varphi \geq 0 \quad \textrm{for all }\,\, x\in M\setminus \partial B_{R_0}\,,
\end{equation}
for some $R_0>0, \lambda>0.$
The non-existence theorem for large initial data also applies for $b=0$. Hence it is a complementary result with respect to those given in \cite{BPT, Pu12, Pu21} for $b=0$.

On the other hand, under suitable assumptions on the vector field $b$, we show that, for any $p>1$, there exists a global solution to problem \eqref{e1}, if the initial datum is small enough. A crucial point in the proof of this result is the construction of a (weak) supersolution to equation
\begin{equation}\label{e36}
\Delta w + \langle b(x), \nabla w \rangle +\lambda w =0 \quad \text{ in }\,\, M\,,
\end{equation}
for suitable $\lambda>0$. 
Observe that differently from the Euclidean case, any $p>1$ is included. Note that such result applies in particular when $b\equiv 0$. In this case, it is in agreement with the results in \cite{Pu12}, \cite{Pu21}.

Consider now the case $(jj)$. We show that for any $p>1$, if $u_0$ is big enough, then problem \eqref{e1} does not admit global solutions. Moreover, we obtain that, for some $\underline p>1$ depending on the Ricci curvature's bound, for every $1<p<\underline p$, for any initial datum $u_0\not\equiv 0$, problem \eqref{e1} cannot have global solutions. In addition, under suitable assumptions on $b$, for some $\bar p>1$ depending on $b$, for every $p>\bar p$, problem \eqref{e1} has global solutions, for initial data $u_0$ small enough.

\medskip

The paper is organized as follows. In section \ref{MB} we recall some useful preliminary notions from Riemannian Geometry and we made the main assumptions. 
Cartan-Hadamard manifolds of type $(j)$ are treated in Section \ref{nonexi}, for non-existence of global solutions, and in Section \ref{exi}, for existence of global solutions.
Moreover, Cartan-Hadamard manifolds of type $(jj)$ are treated in Section \ref{nonexi2}, for non-existence of global solutions, and in Section \ref{exi2}, for existence of global solutions.

\section{Mathematical framework}\label{MB} \setcounter{equation}{0}
\subsection{Preliminaries from Riemannian Geometry}
In this section we recall some useful notions and results from Riemannian Geometry (see e.g. \cite{AMR}, \cite{Grig}, \cite{Grig4}).
%Let $M$ be a complete noncompact Riemannian manifold of dimension $n$. Let $\Delta$
%denote the Laplace-Beltrami operator, $\nabla$ the Riemannian
%gradient and $d\mu$ the Riemannian volume element on $M$.
%In the sequel 
We consider {\it Cartan-Hadamard} manifolds, i.e. simply connected complete noncompact Riemannian manifolds with
nonpositive sectional curvatures. On a Cartan-Hadamard manifold $M$, for any point $o \in M,$
the {\it cut locus } of $o$, $\operatorname{Cut}(o)$, is empty.  Thus $M$ is a manifold with a pole.
Thus, for any $x\in M\setminus \{o\}$ one can define the {\it polar
coordinates} $(r, \theta)$ with respect to $o$, where 
$$r\equiv r(x):=\operatorname{dist}(x,o).$$ %Namely, for
%any point $x\in M\setminus\{o\}$
%there correspond a polar radius $r(x) :=\operatorname{dist}(x, o)$
%and a polar angle $\theta\in \mathbb S^{n-1}$ such that the shortest
%geodesics from $o$ to $x$ starts at $o$ with direction $\theta$ in
%the tangent space $T_o M$. Since we can identify $T_o M$ with
%$\mathbb R^n$, $\theta$ can be regarded as a point of $\mathbb
%S^{n-1}.$ 
For any $x_0\in M$ and for any $R>0$ we set $$B_R(x_0):=
\big\{x\in M\,:\, \operatorname{dist}(x, x_0)<R\,\big\}.$$
We set $B_R\equiv B_R(o).$

The Riemannian metric in $M\setminus\{o\}$ in polar coordinates reads
\[g= dr^2+A_{ij}(r, \theta)d\theta^i d\theta^j, \]
where $(\theta^1, \ldots, \theta^{n-1})$ are coordinates in $\mathbb S^{n-1}$ and $(A_{ij})$ is a positive definite matrix.
%Let $(A^{ij})$ denote the inverse matrix of $(A_{ij})$.
The Laplace-Beltrami operator in polar coordinates has the form
\begin{equation}\label{e6} \Delta = \frac{\partial^2}{\partial r^2} +
\mathcal F(r, \theta)\frac{\partial}{\partial r}+\Delta_{S_{r}},
\end{equation}
where %$\mathcal F(r, \theta):=\frac{\partial}{\partial
%r}\big(\log\sqrt{A(r,\theta)}\big)$, $A(r,\theta):=\det
%(A_{ij}(r,\theta))$, 
\[\mathcal F(r, \theta)=\Delta r(x)\quad \textrm{for any}\;\, x\in M\setminus\{o\}, \]
$\Delta_{S_r}$ is the Laplace-Beltrami operator
on the submanifold $S_{r}:=\partial B_r$\,. The {\it area} element on $S_r$ is 
\[d\mu'|_{S_r}=\sqrt{A} d\theta^1\ldots d\theta^{n-1},\]
where $A:=\operatorname{det}(A_{ij})$, while the {\it volume} element is 
\[d\mu = d\mu' dr\,.\]

A manifold with a pole is a {\it spherically symmetric manifold} or a {\it model}, if the Riemannian metric is given by
\begin{equation}\label{e7}
g= dr^2+\psi^2(r)d\theta^2,
\end{equation}
where %$d\theta^2=\beta_{ij}d\theta^i d \theta^j$ is the standard
%metric in $\mathbb S^{n-1}$, $\beta_{ij}$ being smooth functions of
%$\theta^1, \ldots, \theta^{n-1},$ and 
$\psi\in \mathcal A$ with
$$\mathcal A:=\left\{f\in C^\infty((0,\infty))\cap C^1([0,\infty)): \, f'(0)=1, \, f(0)=0, \, f>0 \ \textrm{in}\;\, (0,\infty)\right\} .$$
%In this case, we write $M\equiv
%M_\psi$; furthermore, we have $\sqrt{A(r,\theta)}=\psi^{n-1}(r)$, so
%that
%\[\Delta = \frac{\partial^2}{\partial r^2}+ (n-1)\frac{\psi'}{\psi}\frac{\partial}{\partial r}+ \frac1{\psi^2}\Delta_{\mathbb S^{n-1}}\,,\]
%where $\Delta_{\mathbb S^{n-1}}$ is the Laplace-Beltrami operator in
%$\mathbb S^{n-1}\,.$ 
Observe that for $\psi(r)=r$, $M=\mathbb R^n$, while for $\psi(r)=\sinh r$, $M$ is the $n-$dimensional hyperbolic space $\mathbb H^n$.

\smallskip

For any $x\in M\setminus\{o\}$, denote by
$\textrm{Ric}_o(x)$ the {\it Ricci curvature} at $x$ in the
radial direction $\frac{\partial}{\partial r}$. 

Let $\omega$ denote any
pair of tangent vectors from $T_xM$ having the form
$\left(\frac{\partial}{\partial r} ,X\right)$, where $X$ is a unit
vector orthogonal to $\frac{\partial}{\partial r}$. We denote by
$\textrm{K}_{\omega}(x)$ the {\it sectional curvature} at the point
$x$ of the plane determined by $\omega$; it is also called {\it sectional radial curvature}.

If $M_\psi$ is a model manifold, then for any $x=(r,
\theta)\in M_\psi\setminus\{o\}$
\[\textrm{K}_{\omega}(x)=-\frac{\psi''(r)}{\psi(r)},\]
and
\[\textrm{Ric}_{o}(x)=-(n-1)\frac{\psi''(r)}{\psi(r)}\,.\]

\subsection{Assumptions}
We always assume that $M$ is a Cartan-Hadamard manifold of dimension $n\ge2$. In addition,  we make one the following assumptions: 

\begin{equation} \tag{{\it $A_0$}} \label{A0}
\begin{cases}
%\textrm{(i)} & M \textrm{ is a Cartan-Hadamard manifold of dimension $n\ge2$} \, ; \\
\textrm{(i)} & \textrm{K}_\omega(x)\leq -h_1^2 \, \textrm{ for some } h_1>0,\, \text{ for any } x\in M\setminus\{o\} ;\\
\textrm{(ii)} &  \textrm{Ric}_o(x)\geq -(n-1)h_2^2 \, \textrm{ for some }\, h_2>h_1,\,  \text{ for any } x\in M\setminus\{o\} ;\\
 
%\;\;\, \textrm{there exists}\;\, \psi\in \mathcal A\,, \hat C>0\,\;\;
%\textrm{such that}\;\frac{\psi'(\rho)}{\psi(\rho)}\leq  \hat C \rho \;\;\textrm{for all}\;\; \rho \ge  1\,, \\
%\qquad \,\, \textrm{Ric}_{o}(x)\geq - (N-1) \frac{\psi''(\rho)}{\psi(\rho)} \quad \textrm{for all}\;\; x\equiv (\rho, \theta)\in M\,.
\end{cases}
\end{equation}
or
\begin{equation} \tag{{\it $A_1$}} \label{A1}
%\begin{cases}
%\textrm{(i)} & M \textrm{ is a Cartan-Hadamard manifold of dimension $n\ge2$} \, ; \\
%\textrm{(ii)} & \textrm{K}_\omega(x)\leq -h_1^2 \text{ for any } x\in M\setminus\{o\},\,\, \textrm{ for some } h_1>0\, ;\\
\textrm{Ric}_o(x)\geq -(n-1)\frac{\bar \beta}{1+ r^2(x)} \, \textrm{ for some }\, \bar\beta>0, \text{ for any } x\in M\setminus\{o\}\,.
 %\;\;\, \textrm{there exists}\;\, \psi\in \mathcal A\,, \hat C>0\,\;\;
%\textrm{such that}\;\frac{\psi'(\rho)}{\psi(\rho)}\leq  \hat C \rho \;\;\textrm{for all}\;\; \rho \ge  1\,, \\
%\qquad \,\, \textrm{Ric}_{o}(x)\geq - (N-1) \frac{\psi''(\rho)}{\psi(\rho)} \quad \textrm{for all}\;\; x\equiv (\rho, \theta)\in M\,.
%\end{cases}
\end{equation}

Concerning the initial datum $u_0$ and the vector field $b$ we always assume:

\begin{equation} \tag{{\it $H$}} \label{H}
\begin{cases}
\textrm{(i) } & b\in C^1(M)\,;\\
\textrm{(ii)} & u_0\in C(M), u_0\geq 0\,.
\end{cases}
\end{equation}

In view of assumption \eqref{A0}, by Laplace and Hessian comparison results (see, e.g., \cite{AMR}, \cite{Grig}),
\begin{equation}\label{e2}
(n-1)h_1\coth(h_1 r) \leq \mathcal F(r, \theta)\leq (n-1)h_2\coth(h_2 r) \quad \textrm{for any } x\equiv (r, \theta)\in M\setminus\{o\}\,.
\end{equation}
Moreover,  there exists $C_0>0$ such that 
\begin{equation}\label{e3}
\mu'(S_R) \leq C_0e^{(n-1) h_2 R}  \quad \textrm{for any}\;\; R>0\,.
\end{equation}

On the other hand, in view of 
 \eqref{A1}, again by Laplace and Hessian comparison results (see, e.g., \cite[Lemma 5.4]{BRS} and remarks at pages 411-412 in \cite{PRS}),
\begin{equation}\label{e2b}
\frac{(n-1)}r\leq \mathcal F(r, \theta)\leq \frac{C_1 (n-1)}{r} \quad \textrm{for any } x\equiv (r, \theta)\in M\setminus\{o\}\,,
\end{equation}
for some $C_1>1$. 
Moreover, there exist $0<C_2<C_3$ 
such that 
\begin{equation}\label{e3b}
C_2  R^{n-1}  \leq \mu'(S_R) \leq C_3 R^{\gamma(n-1)}  \quad \textrm{for any}\;\; R>0\,,
\end{equation}
where 
\begin{equation}\label{e80}
\gamma= \frac{1+\sqrt{1+4\bar\beta}}{2}\,.
\end{equation}

\section{Non-existence of global solutions under assumption \eqref{A0}}\label{nonexi}
Let $a_1>0, C_1>0, a_2>0, C_2>0, R_0>0$. Define 
\[
\tilde \varphi(x)\equiv \tilde \varphi(r):= C_1 e^{- a_1 r} \quad \text{for all }\,\, x\in M\setminus B_{R_0};
\]
\[\hat \varphi(x)\equiv \hat \varphi(r):= C_2 e^{-a_2 r^2}\quad \text{for all }\,\, x\in B_{R_0}\,;\]
\begin{equation}\label{e8}
\varphi:=\begin{cases}
\tilde \varphi & \text{ in  }\,\,  M\setminus B_{R_0}\\
\hat \varphi & \text{ in }\,\, B_{R_0}\,.
\end{cases}
\end{equation}

We always take 
\[C_2=C_1 e^{-a_1 R_0+a_2 R_0^2},\]
\[a_2=\frac{a_1} 2,\]
thus 
\[\varphi\in C^2(M\setminus \partial B_{R_0})\cap C^1(M)\,.\]

As usual, we denote by $\frac{\partial}{\partial r}$ the unitary vector field in the {\it radial} direction. Consider the function $\underline b: M\to \mathbb R$ defined by 
\begin{equation}\label{e50}
\underline b(x):= \left\langle b(x), \frac{\partial}{\partial r}\right\rangle, \quad x\in M\,.
\end{equation}
Clearly, $\underline b(x) \frac{\partial}{\partial r}$ is the {\it radial component} of the vector field $b(x),\, x\in M$.

\begin{lemma}\label{lemma1} 
Assume \eqref{A0}, \eqref{H}. Suppose that there exist $b_0\in \mathbb R, \hat C\geq 0$ such that 
\begin{equation}\label{e10}
\underline{b}(x) \geq  \underline b_0 \quad \text{ for all }\, x\in M\,,
\end{equation}
\begin{equation}\label{e49}
\operatorname{div} b(x)\leq \hat C \quad \text{ for all }\,\, x\in M\,.
\end{equation}
Let $\lambda >\hat C$ and 
\begin{equation}\label{e11}
0<a_1\leq \min\left\{\frac {\lambda-\hat C}{(n-1)h_2\coth(h_2R_0)+\underline{b}_0^-}, \frac{\lambda-\hat C}{1+R_0[(n-1) h_2 \coth(h_2R_0) + \underline{b}_0^-]}\right\}\,. 
\end{equation}
Then the function $\varphi$ defined in \eqref{e8} satisfies \eqref{e4}.
\end{lemma}

\begin{proof}

In view of \eqref{e6}, \eqref{e2} and \eqref{e11} we have, for all $x\in M\setminus \overline{B}_{R_0}$, 
\begin{equation}\label{e12}
\begin{aligned}
&\Delta\varphi -  \langle b(x), \nabla \varphi\rangle -\operatorname{div } b(x) \varphi + \lambda \varphi\\  &=\varphi''(r) + \mathcal F(r, \theta) \varphi'(r) - \underline{b}(x) \varphi'(r) + (\lambda - \hat C) \varphi(r) 
\\& 
\geq e^{-a_1 r}\left\{a_1^2+a_1[-(n-1) h_2 \coth(h_2 R_0) + \underline{b}(x) ] + \lambda-\hat C\right\}\\
& \geq  e^{-a_1 r}\left\{a_1^2+ a_1[-(n-1) h_2 \coth( h_2 R_0) +  \underline{b}_0] + \lambda-\hat C\right\}\\&\geq 
e^{-a_1 r}\left\{-a_1[(n-1) h_2 \coth( h_2 R_0) +  \underline{b}_0^- ] + \lambda-\hat C\right\}
\geq 0\,. 
\end{aligned}
\end{equation}
Furthermore, since $a_2=\frac {a_1}2$, again by \eqref{e6}, \eqref{e2} and \eqref{e11}, for all $x\in B_{R_0}$,
\begin{equation}\label{e13}
\begin{aligned}
&\Delta\varphi - \langle b(x), \nabla \varphi\rangle  -\operatorname{div } b(x) \varphi+ \lambda \varphi\\ &=\varphi''(r) + \mathcal F(r, \theta) \varphi'(r) - \underline b(x) \varphi'(r) +(\lambda-\hat C) \varphi(r) 
\\ &\geq e^{-a_2 r^2}\left\{4 a_2^2 r^2-2a_2+ 2 a_2r [-(n-1) h_2 \coth(h_2 r) + \underline{b}(x) ] +\lambda-\hat C\right\}\\ &
\geq e^{-a_2 r^2} \left\{-2 a_2[1 + R_0((n-1) h_2 \coth(h_2R_0) + \underline{b}_0^-)] +\lambda-\hat C \right\}\geq 0\,.
\end{aligned}
\end{equation}
From \eqref{e12} and \eqref{e13} we get the thesis. 
\end{proof}

Let $a_1>(n-1)h_2,$ so $\varphi\in L^1(M)$. Set
\begin{equation}\label{e41}
c:= \frac1{\|\varphi\|_{L^1(M)}}\,.
\end{equation}

\begin{theorem}\label{thm1}
Assume \eqref{A0}, \eqref{H}, \eqref{e49}.  Let $b\in L^\infty(M)$, $\lambda>\hat C$ and $\varphi$ be as in Lemma \ref{lemma1}
with
\begin{equation}\label{e14}
a_1 >(n-1) h_2\,.
\end{equation}
Let $u$ be a solution of problem \eqref{e1} with $u\in L^\infty\big(M\times (0, \tau)\big)$ for each $0<\tau<T$ and $\partial_t u(t) \, \varphi\in L^1(M)$ for every $t\in (0, T)$. Suppose that 
\begin{equation}\label{e5}
 c \int_M u_0(x) \varphi(x)  d\mu > \lambda^{\frac 1{p-1}}\,,
\end{equation}
with $c$ given by \eqref{e41}. Then problem \eqref{e1} does not admit global solutions. 
\end{theorem}
In Theorem \ref{thm1} the solution to \eqref{e1} is meant in the classical sense. Note that, as it will be apparent from the proof, the same conclusion holds for supersolutions.

\begin{proof}
Consider a family $\{\tilde \zeta_R\}_{R>0}\subset C^\infty_c([0, +\infty))$ of cut-off functions such that, for any $R>0$, 
\[\tilde \zeta_R=\begin{cases}
1 & \text{ in } \left[0, R \right]\\
0 & \text{ in } [2R, +\infty)
\end{cases}\,,\]
\[  0\leq \tilde \zeta_R\leq 1, \quad |\tilde \zeta_R'|\leq \frac{C}{R}, \quad  |\tilde \zeta_R''|\leq  \frac{C}{R^2}  \;\; \text{for any } r\in [0, +\infty)\,,\]
for some $C>0.$
Set
\[\tilde \zeta_R(x) =\tilde \zeta_R(r(x)), \quad x\in M\,.\]
Therefore, 
\begin{equation}\label{e70}
|\nabla \zeta_R(x)|\leq \frac C R\chi_{B_{2R}\setminus B_R}(x)\quad \textrm{ for any } x\in M. 
\end{equation}
Furthermore, in view of \eqref{e6}, \eqref{e2}, 
\begin{equation}\label{e71}
|\Delta \zeta_R(x)| \leq  \frac C R\chi_{B_{2R}\setminus B_R}(x) \quad \textrm{ for any } x\in M. 
\end{equation}

From \eqref{e3} and \eqref{e14} we can infer that $\varphi\in L^1(M)$. 
Define 
\[\psi(t):=\int_M c\,  u(x,t) \varphi(x)\, d\mu, \quad t\in (0, T)\,. \]

From \eqref{e1} we obtain 
\begin{equation}\label{e16}
\int_M c \varphi\zeta_R \pa_t u d\mu = \int_M c \varphi \zeta_R \Delta u d\mu + \int_M c \varphi\zeta_R \langle b, \nabla u\rangle d\mu + \int_M c\zeta_R \varphi u^p d\mu\,. 
\end{equation}
Let $\nu$ be the outer normal unit vector to $\pa B_{R_0}$. Since $\varphi \in C^2(M\setminus\pa B_{R_0})\cap C^1(M)$, we get
\begin{equation}\label{e18}
\begin{aligned}
\int_M \varphi \zeta_R \Delta u d\mu&=\int_{B_{R_0}}\varphi \zeta_R \Delta u d\mu + \int_{B_{R_0}^c}\varphi \zeta_R \Delta u d\mu \\ 
&= -\int_{B_{R_0}}\langle \nabla u, \nabla \varphi_1 \rangle \zeta_R d\mu +\int_{\pa B_{R_0}}\varphi_1 \zeta_R \frac{\partial u}{\partial \nu} dS\\
&-\int_{B^c_{R_0}}   \langle \nabla u, \nabla \varphi_2 \rangle \zeta_R d\mu - \int_{\pa B_{R_0}}\varphi_2 \zeta_R \frac{\partial u}{\partial \nu} dS\\& -\int_M \langle \nabla u, \nabla \zeta_R \rangle \varphi d\mu \\
&=\int_{B_{R_0}}u \Delta\varphi_1 \zeta_R d\mu -\int_{\partial B_{R_0}} u\zeta_R \frac{\partial \varphi_1}{\pa \nu} dS \\
& +\int_{B_{R_0}^c} u \Delta\varphi_2 \zeta_R d\mu +\int_{\partial B_{R_0}} u\zeta_R \frac{\partial \varphi_2}{\pa \nu} dS\\
& + 2 \int_M u \langle \nabla\varphi, \nabla \zeta_R\rangle  d\mu +\int_M u \varphi \Delta\zeta_R d\mu\\
&=\int_{B_{R_0}}u \Delta\varphi_1 \zeta_R d\mu   +\int_{B_{R_0}^c} u \Delta\varphi_2 \zeta_R d\mu \\
& + 2 \int_M u \langle \nabla\varphi, \nabla \zeta_R\rangle  d\mu +\int_M u \varphi \Delta\zeta_R d\mu\\
&=  \int_{B_{R_0}} u \zeta_R\Delta\varphi_1 d\mu +   \int_{B^c_{R_0}} u \zeta_R \Delta\varphi_2 d\mu \\
& + 2 \int_M u \langle \nabla\varphi, \nabla \zeta_R\rangle  d\mu +\int_M u \varphi \Delta\zeta_R d\mu.\\
\end{aligned}
\end{equation}
Observe that in view of our assumptions, Lemma \ref{lemma1} can be applied. Since $u\geq 0, \varphi>0, \zeta_R\geq 0$, from \eqref{e4} and \eqref{e18} we obtain
\begin{equation}\label{e19}
\begin{aligned}
\int_M \varphi \zeta_R \Delta u d\mu
&\geq -\lambda \int_{B_{R_0}} u \varphi_1 \zeta_R d\mu + \int_{B_{R_0}}\langle b(x), \nabla \varphi_1\rangle u \zeta_R d\mu \\
&-\lambda \int_{B^c_{R_0}} u \varphi_2 \zeta_R d\mu + \int_{B^c_{R_0}}\langle b(x), \nabla \varphi_2\rangle u \zeta_Rd\mu\\
&+ \int_M u \varphi \zeta_R \operatorname{div} b\, d\mu+ 2 \int_M u \langle \nabla\varphi, \nabla \zeta_R\rangle  d\mu +\int_M u \varphi \Delta\zeta_R d\mu\\
&=  -\lambda \int_{M} u \varphi \zeta_R d\mu + \int_{M}\langle b(x), \nabla \varphi\rangle u \zeta_Rd\mu+\int_M u \varphi \zeta_R \operatorname{div} b\, d\mu\\
& + 2 \int_M u \langle \nabla\varphi, \nabla \zeta_R\rangle  d\mu +\int_M u \varphi \Delta\zeta_R d\mu\,.
\end{aligned}
\end{equation}
Furthermore, 
\begin{equation}\label{e17}
\begin{aligned}
&\int_M \varphi \zeta_R \langle b, \nabla u\rangle d\mu\\ & = -\int_M u \varphi \zeta_R \operatorname{div} b\, d\mu - \int_M u \zeta_R \langle b, \nabla\varphi \rangle d\mu 
-\int_M u \varphi \langle b, \nabla\zeta_R \rangle d\mu.
\end{aligned}
\end{equation}
By \eqref{e16}, \eqref{e19} and \eqref{e17},
\begin{equation*}
\begin{aligned}
\int_M c \varphi\zeta_R \pa_t u d\mu& \geq - \int_{M} \lambda u \varphi \zeta_R d\mu +  \int_{M}\langle b(x), \nabla \varphi\rangle  u \zeta_Rd\mu\\
& + 2 \int_M  u \langle \nabla\varphi, \nabla \zeta_R\rangle  d\mu +\int_M u \varphi \Delta\zeta_R d\mu\\
&  - \int_M  \langle b(x), \nabla\varphi \rangle u \zeta_R d\mu 
-\int_M\langle b(x), \nabla\zeta_R \rangle  u \varphi  d\mu\\ & + \int_M c\zeta_R \varphi u^p d\mu\,.
\end{aligned}
\end{equation*} 
Therefore,
\begin{equation}\label{e20}
\begin{aligned}
\int_M c \varphi\zeta_R \pa_t u d\mu & +\lambda \int_{M} u \varphi \zeta_R d\mu  \geq  
 2 \int_M u \langle \nabla\varphi, \nabla \zeta_R\rangle  d\mu\\ & +\int_M u \varphi \Delta\zeta_R d\mu
-\int_M u \varphi \langle b, \nabla\zeta_R \rangle d\mu + \int_M c\zeta_R \varphi u^p d\mu\,.
\end{aligned}
\end{equation} 
In view of \eqref{e3}, \eqref{e70} and \eqref{e14}, for every $t\in (0, T)$, we get, for some $C>0$, 
\begin{equation}\label{e21}
\begin{aligned}
\left|\int_M  u(x,t) \langle \nabla\varphi, \nabla \zeta_R\rangle  d\mu\right| & \leq \frac{\|u(t)\|_{\infty} C}{R}\int_{B_{2R}} |\nabla \varphi|d\mu \\ &\leq \frac{\|u(t)\|_{\infty} C}{R}\int_0^{2R} e^{[(n-1)h_2-a_1]r}dr \mathop{\longrightarrow}_{R\to+\infty} 0\,.
\end{aligned}
\end{equation}
From \eqref{e3}, \eqref{e71} and \eqref{e14},  for every $t\in (0, T)$, we obtain, for some $C>0$, 
\begin{equation}\label{e22}
\begin{aligned}
\left|\int_M u \varphi \Delta\zeta_R d\mu\right| \leq \frac{\|u(t)\|_{\infty} C}{R}\int_0^{2R} e^{[-a_1+(n-1)h_2]r }dr \mathop{\longrightarrow}_{R\to+\infty} 0\,.
\end{aligned}
\end{equation}
In view of \eqref{e3} and \eqref{e14}, for every $t\in (0, T)$, we get, for some $C>0$, 
\begin{equation}\label{e23}
\left|\int_M u \varphi \langle b, \nabla\zeta_R \rangle d\mu\right| \leq 
\frac{\|u(t)\|_{\infty} \|b\|_{\infty}C}{R}\int_0^R e^{[(n-1)h_2-a_1]r}dr \mathop{\longrightarrow}_{R\to+\infty} 0\,.
\end{equation}
Since, for every $t\in (0, T)$, $u(t)\in L^\infty(M)$ and $\varphi\in L^1(M)$, by the dominated convergence theorem, 
\begin{equation}\label{e24b}
\int_{M} u \varphi \zeta_R d\mu  \mathop{\longrightarrow}_{R\to+\infty} \int_M c \varphi u\, d\mu=\psi(t)\,,
\end{equation}
and
\begin{equation}\label{e24}
 \int_M c\zeta_R \varphi u^p d\mu \mathop{\longrightarrow}_{R\to+\infty} \int_M c \varphi u^p d\mu\geq \psi^p(t)\,;
\end{equation}
here use of Jensen's inequality and \eqref{e41} have been made. 
Since, for every $t\in (0, T)$, $\partial_t u(t) \varphi \in L^\infty(M)$, by the dominated convergence theorem, 
\begin{equation}\label{e25}
\int_M c \varphi\zeta_R \pa_t u\, d\mu \mathop{\longrightarrow}_{R\to+\infty} \int_M c\varphi \pa_t u\, d\mu = \psi'(t)\,.
\end{equation}
By \eqref{e20}-\eqref{e25} and \eqref{e5}
\[\psi'(t) \geq \psi^p(t)- \lambda \psi(t)  \,,\quad \psi(0)>  \lambda^{\frac 1{p-1}}\,.
\]
This easily implies that there exists $T>0$ such that 
\[\lim_{t\to T^-} \psi(t)=+\infty\,.\]
Hence the thesis follows.
\end{proof}

\section{Existence of global solutions under assumption \eqref{A0}}\label{exi} \setcounter{equation}{0}
Let $\underline b(x)$ be defined as in \eqref{e50}. Assume that 
\begin{equation}\label{e37}
\underline b(x)\geq \underline b_1> - (n-1) h_1 \quad \text{ for all } \, x\in M\,. 
\end{equation}
Set $$\sigma:=(n-1)h_1+\underline b_1.$$ 
So, $\sigma>0$.
Suppose that 
\begin{equation}\label{e38}
0<\lambda < \frac{\sigma^2}{4}\,.
\end{equation}
Hence we can find $a\in \mathbb R$ such that 
\begin{equation}\label{e39}
0<a<\frac{\sigma+\sqrt{\sigma^2-4\lambda}}{2}\,.
\end{equation}

Define 
\[w(x)\equiv w(r):= e^{-a r}, \quad x\in M\,.\]

\begin{proposition}\label{prop1}
Let assumptions \eqref{A0} and \eqref{H} be satisfied. Let conditions \eqref{e37}-\eqref{e39} be fulfilled.
Then $w$ is a weak supersolution of \eqref{e36}.
\end{proposition}
\begin{proof}
In view of \eqref{e2} and \eqref{e39}, for any $x\in M\setminus\{o\}$,
\begin{equation}\label{e40}
\begin{aligned}
&\Delta w + \langle b(x), \nabla w(x) \rangle + \lambda w(x)\\ & = e^{- ar}\left[a^2- a\mathcal F(r, \theta) - a \underline b(x) +\lambda \right]\\
& \leq e^{-a r}\left\{a^2-a[(n-1)h_1+\underline b_1]+\lambda \right\}\leq 0\,.
\end{aligned}
\end{equation}
Since $w'(0)<0,$ by a Kato's type inequality, $w$ is a weak supersolution of equation \eqref{e36} in the whole $M$. 
\end{proof}

\begin{theorem}\label{thm2}
Let assumptions \eqref{A0}, \eqref{H}, \eqref{e37} and \eqref{e38} be satisfied. Let $w$ be defined as in Proposition \ref{prop1}. 
Assume that
\begin{equation}\label{e30}
0\le u_0\le \tilde C w \quad \textrm{in} \;\; M\,,
\end{equation}
where
\begin{equation}\label{e31}
0< \tilde C< \frac 1{\|w \|_{\infty}} \lambda^{\frac 1{p-1}}\,.
\end{equation}
Then there exists a global solution $u$ of problem
\eqref{e1}; in addition, $u\in L^\infty(M\times(0, +\infty))\,.$
\end{theorem}

\bigskip
\bigskip

In order to prove Theorem \ref{thm2}, we adapt to the present situation some arguments used in \cite[Theorem 3.1]{BPT} (see also \cite[Theorem 3.2]{Pu21})\,.  Let $\{\Omega_j\}_{j\in\ene}$ be a sequence of domains
$\{\Omega_j\}_{j\in \ene}\subseteq M$ such that
$\bar\Omega_j\subseteq \Omega_{j+1}$ for every $j\in \ene,\,
\bigcup_{j=1}^\infty\Omega_j=M\,,\pa \Omega_j$ is smooth for every
$j\in \ene\,.$ Furthermore, for every $j\in \mathbb N$ let $\zeta_j\in C^\infty_c(\Omega_j)$ such that $0\leq \zeta_j\leq 1, \zeta_j\equiv 1$ in $\Omega_{j/2}.$

\begin{proof}[Proof of Theorem \ref{thm2}\,.]  For any  $j\in \ene$ there esists a unique classical solution $u_j$ to
problem
\begin{equation}\label{e33}
\left\{
\begin{array}{ll}
 \,  \pa_t u = \Delta u \,+\, \langle b(x), \nabla u \rangle +  u^p &\textrm{in}\,\,\Omega_j\times (0,T)
\\&\\
\textrm{ }u \,= 0 & \textrm{in}\,\, \pa\Omega_j\times (0,T) \\&\\
\textrm{ }u \, = \zeta_j u_0& \textrm{in}\,\,  \Omega_j\times \{0\} \,.
\end{array}
\right.
\end{equation}

\smallskip
Take the constant $\tilde C>0$ given by \eqref{e31}. Let
\[\tilde w(x):= \tilde C w(x)\quad \big(x\in M\big)\,,\]
and
\[\xi(t)=\left\{1-\frac 1{\lambda}\|\tilde w \|_{\infty}^{p-1}\left[ 1-e^{-(p-1)\lambda t}\right] \right\}^{-\frac 1{p-1}}\,\quad \big(t\in [0,\infty)\big)\,.\]

Note that $\xi$ is well-defined in $[0, \infty)$ due to \eqref{e31}.
It is easily seen that $\xi$ solves
problem
\begin{equation}\label{e41a}
\left\{
\begin{array}{ll}
 \,  \xi' = \|\tilde w \|_\infty^{p-1}e^{-(p-1) \lambda t} \xi^p\,,  &
 t\in (0,\infty)
\\&\\
\textrm{ }\xi(0) \, = 1 \,.
\end{array}
\right.
\end{equation}
Select any $\lambda$ satisfying \eqref{e38}. Define
\begin{equation}\label{e41b}
\bar u(x,t):= e^{-\l t}\xi(t)\tilde w (x)\quad \big((x,t)\in M\times[0,\infty)\big)\,.
 \end{equation}
Due to \eqref{e41a}, we have 
\[\partial_t \bar u - \Delta \bar u - \langle b(x), \nabla \bar u\rangle - \bar u^p \]\[= -\l e^{-\l t}\xi(t) \tilde w(x) + e^{-\l t}\|\tilde w \|_{\infty}^{p-1} e^{-(p-1)\l t} \xi^p(t) \tilde w(x) \]\[+ \l e^{-\l t}\xi(t) \tilde w(x) -e^{-\l p t}\xi^p(t) \tilde w^p (x) \geq 0\quad  \textrm{weakly in}\,\, M\times (0, \infty)\,. \]
So, $\bar u$ is a weak supersolution of equation
\begin{equation}\label{e42}
\pa_t u = \Delta  u + \langle b(x), \nabla u \rangle+ u^p \quad \textrm{in} \;\; M\times (0, \infty)\,.
\end{equation}
Moreover, due to \eqref{e30}, for any $j\in \ene$, $\bar u$ is a bounded weak
supersolution of problem \eqref{e33}. Obviously, for each $j\in \ene$, $\underline u\equiv 0$ is a subsolution of problem \eqref{e33}.
Hence, by the comparison principle, for every $j\in \ene$ we obtain
\begin{equation}\label{e43}
0 \le u_j \le \bar u\quad \textrm{in}\;\; B_j\times (0, T)\,.
\end{equation}
%Fix any $N\in \mathbb N.$ For all $j\in \mathbb N, j\geq N+1$, for every $\varepsilon \in (0, T)$, in view of \eqref{e29} and Remark \ref{partheory}-(b) with $u=u_j, K=B_N$, we can infer that
%\[\begin{aligned}
%&\|D^2u_j \|_{C_{x,t}^{\sigma, \sigma/2}(B_N\times [\varepsilon, T])}+ \|D u_j \|_{C_{x,t}^{\sigma, \sigma/2}(B_N\times [\varepsilon, T])}  \\ &+\|\pa_t u_j\|_{C_{xt}^{\sigma, \sigma/2}(B_N\times [\varepsilon, T])}+ \|u_j\|_{C^{\sigma, \sigma/2}_{x,t}(B_N\times [0, T])}\leq C_{N, \varepsilon}\,,
%\end{aligned}
%\]
%where $C_{N, \varepsilon}$ depends on $\|\bar u\|_{L^{\infty}(B_N\times (0, T))}$, but it is independent of $j$. By Ascoli-Arzel\`a Theorem, we can extract a subsequence $\{u_{j_k}\}$ from $\{u_j\}$, which converges in $C_{x,t}^{2,1}(M\times (0, T])$ and in $C(M\times [0, T])$. Since we can repeat the argument for every $N$, we can perform a diagonal argument to construct a subsequence $\{u_{j_{k_h}}\}$ of $\{u_{j_k}\}$, which converges in $C^{2,1}_{x,t}(K\times [\varepsilon, T])$ as $h\to +\infty$, for each compact subset $K\subset M$ and for each $\varepsilon\in (0, T)$, and in $C_{\rm loc}(M\times [0, T])$, to some function $u\in C_{x,t}^{2,1}(M\times (0, T])\cap C(M\times [0, T])$. Moreover, $u$ is a classical solution problem \eqref{e1}.
By standard a priori estimates (see, e.g., \cite{Fr}, \cite{LSU}), we can infer that there exists a subsequence $\{u_{j_{k_h}}\}$ of $\{u_{j_k}\}$, which converges in $C^{2,1}_{x,t}(K\times [\varepsilon, T])$ as $h\to +\infty$, for each compact subset $K\subset M$ and for each $\varepsilon\in (0, T)$, and in $C_{\rm loc}(M\times [0, T])$, to some function $u\in C_{x,t}^{2,1}(M\times (0, T])\cap C(M\times [0, T])$. Moreover, $u$ is a classical solution problem \eqref{e1}.
Furthermore, from \eqref{e43} we get
\[0\le u\le \bar u\quad\textrm{in}\;\;M\times (0, T)\,.
\]
Hence the thesis follows.
\end{proof}

%Clearly, if $h$ and $u_0$ are H\"older continuous, then, in view of standard methods (see e.g. \cite{Fr}, \cite{LSU}, \cite{L}), we have that $u_{j_k}$ converges locally in $C^{2,1}_{x,t}(M\times [0, T)$ to a solution of problem \eqref{e1}.
%
%\smallskip
%
%In general the solution of problem \eqref{e1} is not unique in $L^{\infty}(M\times (0, T))$. Conditions related to $M$ that guarantee uniqueness for problem \eqref{e1} and the validity of comparison principles are established e.g. in \cite[Corollary 3.7]{Pu12} (see also e.g. \cite{Grig} for linear equations).
%
%\begin{remark}\label{oss3} In \cite{Pu12} in the proof of global existence it is used a comparison principle in $M\times (0, T)$; this holds under suitable assumptions on $M$.
%Indeed, here we use comparison principles only in $B_j\times (0, T)$ for every $j\in \ene$, hence we do not require those assumptions. Moreover, observe that in particular in Theorem \ref{thm4} when \eqref{A0} is satisfied with $\g>2$, the comparison principle on $M\times (0, T)$ does not hold.
%\end{remark}

\section{Non-existence of global solutions under assumption \eqref{A1}}\label{nonexi2}
Let $a>0$. Define 
\begin{equation}\label{e53}
\eta(x)\equiv \eta(r):= e^{- a r^2} \quad \text{for all }\,\, x\in M.
\end{equation}

\begin{lemma}\label{lemma2} 
Assume \eqref{A1}, \eqref{H}. Suppose that \eqref{e49} holds and that there exists $\sigma>0$ such that 
\begin{equation}\label{e51}
\underline{b}(x) \geq  - \frac{\sigma}{r(x)} \quad \text{ for all }\, x\in M\setminus\{o\}\,,
\end{equation}
with $\underline b(x)$ defined in \eqref{e50}. Let $a>0$ and 
\begin{equation}\label{e52}
\lambda \geq 2a [1+ C_1(n-1)+\sigma] + \hat C\,,
\end{equation}
with $C_1$ given by \eqref{e2b}. Then the function $\eta$ defined in \eqref{e53} satisfies \eqref{e4}.
\end{lemma}

\begin{proof}

In view of \eqref{e6}, \eqref{e2b}, \eqref{e51} and \eqref{e52} we have, for all $x\in M\setminus \overline{B}_{R_0}$, 
\begin{equation*}
\begin{aligned}
&\Delta\eta -  \langle b(x), \nabla \eta\rangle -\operatorname{div } b(x) \eta + \lambda \eta\\  &\geq \eta''(r) + \mathcal F(r, \theta) \eta'(r) - \underline{b}(x) \eta'(r) + (\lambda - \hat C) \varphi(r) 
\\& 
\geq e^{-a r^2}\left\{ 4a^2 r^2 - 2a  - 2a r \frac{(n-1) C_1}{r} + 2a r \underline{b}(x)  + \lambda-\hat C\right\}\\
& \geq 
e^{-a r^2}\left\{-2a[1+ C_1(n-1) +\sigma ] + \lambda-\hat C\right\}
\geq 0\,. 
\end{aligned}
\end{equation*}
\end{proof}

Due to \eqref{e3b}, $\eta\in L^1(M)$. Set
\begin{equation}\label{e55}
k:= \frac1{\int_M \eta(x) d\mu}\,.
\end{equation}

\begin{theorem}\label{thm3}
Assume \eqref{A1}, \eqref{H}, \eqref{e49}, \eqref{e51}, \eqref{e52}.  Let $b\in L^\infty(M)$  and $\varphi$ be as in Lemma \ref{lemma2}.
Let $u$ be a solution of problem \eqref{e1} with $u\in L^\infty\big(M\times (0, \tau)\big)$ for each $0<\tau<T$ and $\partial_t u(t) \, \varphi\in L^1(M)$ for every $t\in (0, T)$. Suppose that 
\begin{equation}\label{e56}
 k \int_Mu_0(x)  \eta(x)  d\mu > \lambda^{\frac 1{p-1}}\,,
\end{equation}
with $k$ given by \eqref{e55}. Then problem \eqref{e1} does not admit global solutions. 
\end{theorem}
The proof of Theorem \ref{thm3} is very similar to that of Theorem \ref{thm1}, and it is based on Lemma \ref{lemma2} instead of Lemma \ref{lemma1}. Therefore it is omitted. 

\begin{corollary}\label{cor1}
Assume \eqref{A1}, \eqref{H}, and \eqref{e49}, with $\hat C=0$.  Let $u_0\not\equiv 0, b\in L^\infty(M)$. 
Let $u$ be a solution of problem \eqref{e1} with $u\in L^\infty\big(M\times (0, \tau)\big)$ for each $0<\tau<T$ and $\partial_t u(t) \, \varphi\in L^1(M)$ for every $t\in (0, T)$.
Suppose that 
\begin{equation}\label{e57}
1<p<1+\frac 2{\gamma(n-1)+1}\,,
\end{equation}
with $\gamma$ given by \eqref{e80}. Then problem \eqref{e1} does not admit global solutions. 
\end{corollary}

\begin{remark}
Observe that if $M=\mathbb R^n$, then \eqref{A1} is fulfilled with $\bar \beta=0$. So, $\gamma=1$. Therefore, condition \eqref{e57} gives 
\[1<p<1+\frac 2 n\,.\]
This agrees with the  results in \cite{Fuji}, \cite{Hay}, \cite{KST}. Indeed, in $\mathbb R^n$ also the equality sign included. However, by means of our methods we are not able to consider the equality sign in condition \eqref{e57}.
\end{remark}

\begin{proof} The conclusion follows from Theorem \ref{thm3}, if we show that condition \eqref{e56} is fulfilled for $a>0$ and 
\[\lambda=2a [1+ C_1(n-1)+\sigma]\,. \]

Note that, for some $\tilde C_1>0$ (see, e.g., \cite[formula (3.3.11)]{BL}), for any $a>0$,
\[\int_0^\infty e^{-a r^2} r^{\gamma(n-1)} dr \leq  \tilde C_1 a^{-\frac n2-\frac{(\gamma-1)(n-1)}{2}}=\tilde C_1 a^{-\frac{\gamma(n-1)+1}2}\,. \]
Hence, in view of \eqref{e3b} and \eqref{e55}, for some $\tilde C_2>0$, 
\[1=k \int_M \eta(x) d\mu\leq k \tilde C_2 \int_0^\infty e^{-a r^2} r^{\gamma(n-1)} dr \leq k \tilde C_2 \tilde C_1 a^{-\frac{\gamma(n-1)+1}2}\,.  \]
Thus
\[k\geq \frac{a^{\frac{\gamma(n-1)+1}2}}{\tilde C_1\tilde C_2}\,,\]
and
\begin{equation}\label{e59}
k \int_M \eta(x) u_0(x) d\mu \geq \frac{a^{\frac{\gamma(n-1)+1}2}}{\tilde C_1\tilde C_2}\int_M \eta(x) u_0(x) d\mu \,.
\end{equation}
Hence condition \eqref{e56} is verified, if 
\[ \frac{a^{\frac{\gamma(n-1)+1}2}}{\tilde C_1\tilde C_2}\int_M \eta(x) u_0(x) d\mu > a^{\frac 1{p-1}} [2(1+ C_1(n-1)+\sigma)]^{\frac 1{p-1}}\,,\]
that is
\begin{equation}\label{e60}
\int_M \eta(x) u_0(x) d\mu >\tilde C_1\tilde C_2\, a^{-\frac{\gamma(n-1)+1}2+\frac 1{p-1} }    [2(1+ C_1(n-1)+\sigma)]^{\frac 1{p-1}}\,.
\end{equation}
Since $u_0\not\equiv 0$, 
\[\liminf_{a\to 0^+}\int_M \eta(x) u_0(x) d\mu  >0; \]
furthermore, in view of \eqref{e57}, we have  
\[
\lim_{a\to 0^+} a^{-\frac{\gamma(n-1)+1}2+\frac 1{p-1} }[2(1+ C_1(n-1)+\sigma)]^{\frac 1{p-1}}=0\,.
\]
Hence, conditions \eqref{e60} and \eqref{e56}  are satisfied for $a>0$ sufficiently small. This completes the proof. 
\end{proof}

\section{Existence of global solutions under assumption \eqref{A1}}\label{exi2}
Let $\underline b(x)$ be defined as in \eqref{e50}. For any $C>0, \alpha>0, t_0>0$ define 
\begin{equation}\label{e71}
\underline u(x,t):=C(t+t_0)^{-\alpha}e^{-\frac{r^2}{4(t+t_0)}} \quad \text{ for any }\, x\in M, t>0\,.
\end{equation}

\begin{proposition}\label{prop2}
Let assumptions \eqref{A1} and \eqref{H} be satisfied. Suppose that 
\begin{equation}\label{e70}
\underline b(x) \geq \frac{\nu}{r(x)}\quad \text{ for all }\,\, x\in M\setminus\{o\},
\end{equation}
for some $-n<\nu\leq 0$, and that 
\begin{equation}\label{e72}
p>1+\frac{2}{n+\nu}\,.
\end{equation}
 Then, for some $\alpha>0, C>0$, the function $\bar u$, defined in \eqref{e71}, verifies
\[\pa_t \bar u \geq  \Delta \bar u + \langle b, \nabla \bar u\rangle + \bar u^p\quad \text{ in }\,\, M\times (0, +\infty)\,.\]   
\end{proposition}
The proof of Proposition \ref{prop2} is modelled after that of \cite[Lemma 3.3.2]{BL}.
\begin{proof}
From \eqref{e6} and \eqref{e2b} we get, for every $(x,t)\in M\times (0, +\infty),$ 
\begin{equation}\label{e73}
\begin{aligned}
&\pa_t \bar u - \Delta \bar u - \langle b(x), \nabla \bar u\rangle  \\ &=C e^{-\frac{r^2}{4(t+t_0)}} (t+t_0)^{-\alpha-1}\left\{-\alpha +\frac1 2+ \frac{r}{2}\mathcal F(r, \theta)+\frac 1 2 \underline b(x) r \right\}\\
&=C e^{-\frac{r^2}{4(t+t_0)}} (t+t_0)^{-\alpha-1}\left\{-\alpha +\frac1 2+ \frac{r}{2}\frac{n-1}{r}+\frac{\nu}2 \right\}\\
& \geq C e^{-\frac{r^2}{4(t+t_0)}} (t+t_0)^{-\alpha-1}\frac{n+\nu-2\alpha}{2}=\frac{\bar u}{t+t_0}\frac{n+\nu-2\alpha}{2}\,.
\end{aligned}
\end{equation}
Due to \eqref{e72} we can find 
\[0<\epsilon<\frac{n+\nu}{2}- \frac1{p-1}\,.\]
Fix any $\epsilon\in \left(0,\frac{n+\nu}{2}\right)$ and choose
\begin{equation}\label{e75}
\alpha=\frac{n+\nu}{2}-\epsilon\,.
\end{equation}
Note that 
\[\bar u(x,t)\leq C (t+t_0)^{-\alpha}\quad \text{ for all }\, (x,t)\in M\times (0, +\infty)\,.\]
Therefore
\begin{equation}\label{e74}
\frac1{t+t_0}\geq \left(\frac{\bar u(x,t)}C \right)^{\frac 1{\alpha}}\quad \text{ for all }\, (x,t)\in M\times (0, +\infty)\,.
\end{equation}
Choose 
\begin{equation}\label{e76}
0< C\leq \epsilon^{\frac 1{\alpha}}\,.
\end{equation}
From \eqref{e73}-\eqref{e76} we obtain 
\begin{equation}\label{e77}
\begin{aligned}
&\pa_t \bar u - \Delta \bar u - \langle b, \nabla \bar u\rangle  \geq \bar u^{1+\frac1{\alpha}}\quad \text{ in }\,\, M\times (0, +\infty)\,.
\end{aligned}
\end{equation}
In view of \eqref{e72} and \eqref{e75}, for $\epsilon>0$ sufficiently small
\begin{equation}\label{e78}
1+\frac1{\alpha}< p\,;
\end{equation}
in addition,
\begin{equation}\label{e79}
0<\bar u \leq 1 \quad \textrm{ in }\,\,  M\times (0, +\infty)\,.
\end{equation}
By \eqref{e77}-\eqref{e79},
\[\pa_t \bar u - \Delta \bar u - \langle b, \nabla \bar u\rangle\geq \bar u^p \quad \textrm{ in }\,\, M\times (0, +\infty)\,.\]
The proof is complete. 
\end{proof}

\begin{theorem}\label{thm4}
Let assumptions \eqref{A1} and \eqref{H} be satisfied. Suppose that \eqref{e70} and \eqref{e72} hold. Assume that 
\[0\leq u_0(x)\leq \bar u(x, 0) \quad \text{ for all }\,\, x\in M,\]
where $\bar u$ given by Proposition \ref{prop2}. Then there exists a global solution to problem \eqref{e1}. 
\end{theorem}
\begin{proof}
The conclusion follows arguing as in the proof of Theorem \ref{thm2} replacing the supersolution defined in \eqref{e41} with that provided by Proposition \ref{prop2}.
\end{proof}

Theorem \ref{thm4} generalizes \cite[Theorem 3.3.3]{BL}, where $M=\mathbb R^n$ is treated.

\end{document}